\def\Time-stamp: <#1 h>{#1}
\documentclass[11pt]{amsart}
\usepackage{amscd}
\usepackage{amssymb}
\theoremstyle{plain}
\newtheorem{theorem}{Theorem}[section]
\newtheorem{lemma}[theorem]{Lemma}
\newtheorem{corollary}[theorem]{Corollary}
\newtheorem{proposition}[theorem]{Proposition}
\newtheorem{conjecture}[theorem]{Conjecture}
\theoremstyle{definition}
\theoremstyle{remark}
\newtheorem{remark}[theorem]{Remark}

\newcommand\modZ {{\mathbb Z}}
\newcommand\modQ {{\mathbb Q}}
\newcommand\modN {{\mathbb N}}
\newcommand\modM {{\mathcal M}}
\newcommand\modF {{\mathcal F}}
\newcommand\PHI{{\mathbf{\Phi }}}

\newcommand\simeqto{{\overset{\simeq}{\longrightarrow}}}

\newcommand\ord{\operatorname{ord}}
\newcommand\id{\operatorname{id}}

\newcommand\zzzvert {\;|\;}
\newcommand\zzzcolon {\colon\thinspace}

\newcommand\xq[2]{{#1[q]^{#2}}}
\newcommand\Rq[1]{\xq{R}{#1}}
\newcommand\Aq[1]{\xq{A}{#1}}
\newcommand\Qq[1]{\xq{\modQ }{#1}}
\newcommand\Zq[1]{\xq{\modZ }{#1}}
\newcommand\RqS{\Rq{S}}
\newcommand\AqS{\Aq{S}}
\newcommand\ZqS{\Zq{S}}
\newcommand\QqS{\Qq{S}}
\newcommand\RqZ{\Rq{Z}}
\newcommand\Rqn{\Rq{\{n\}}}
\newcommand\Qqn{\Qq{\{n\}}}

\newcommand\ZqN{\Zq{\modN }}
\newcommand\RqN{\Rq{\modN }}
\newcommand\QqN{\Qq{\modN }}
\newcommand\RqM{\Rq{M}}

\begin{document}

\title{Cyclotomic completions of polynomial rings}
\author{Kazuo Habiro}
\address{Research Institute for Mathematical Sciences, Kyoto
  University, Kyoto, 606-8502, Japan}
\email{habiro@kurims.kyoto-u.ac.jp}

\begin{abstract}
  The main object of study in this paper is the completion
  $\ZqN=\varprojlim_{n}\modZ [q]/((1-q)(1-q^2)\cdots(1-q^n))$ of the
  polynomial ring $\modZ [q]$, which arises from the study of a new
  invariant of integral homology $3$-spheres with values in $\ZqN$
  announced by the author, which unifies all the $sl_2$
  Witten-Reshetikhin-Turaev invariants at various roots of unity.  We
  show that any element of $\ZqN$ is uniquely determined by its
  power series expansion in $q-\zeta $ for each root $\zeta $ of unity.  We
  also show that any element of $\ZqN$ is uniquely determined by its
  values at the roots of unity.  These results may be interpreted that
  $\ZqN$ behaves like a ring of ``holomorphic functions defined on the
  set of the roots of unity''.  We will also study the generalizations
  of $\ZqN$, which are completions of the polynomial ring $R[q]$ over
  a commutative ring $R$ with unit with respect to the linear
  topologies defined by the principal ideals generated by products of
  powers of cyclotomic polynomials.
\end{abstract}



\date{September 6, 2002}

\maketitle

\tableofcontents

\section{Introduction}
\label{sec:introduction}
The main object of study in this paper is the completion
\begin{equation}
  \label{eq:3}
  \ZqN = \varprojlim_{n\ge 0} \modZ [q] /((q)_n)
\end{equation}
of the polynomial ring $\modZ [q]$ in an indeterminate $q$, where we use
the notation
\[
(q)_n =(1-q)(1-q^2)\cdots(1-q^n)\in \modZ [q], \quad n\ge 0.
\]
Each element $a\in \ZqN$ is expressed as an infinite sum
\begin{equation}
  \label{eq:17}
  a = \sum_{n\ge 0} a_n(q)_n,
\end{equation}
where $a_n\in \modZ [q]$ for $n\ge 0$.

Some specific instances of the series~\eqref{eq:17} and some variants
which can define elements of $\ZqN$ can be found in the literature.
Zagier \cite{Zagier} studied the series $\sum_{n\ge 0}(q)_n$, which was
introduced by Kontsevich, and observed that it can be expanded in
$q-\zeta $ for any root $\zeta $ of unity.  Clearly, this is the case also for
any elements of $\ZqN$.  Some of the formulae given by Lawrence and
Zagier~\cite{Lawrence-Zagier} and by Le~\cite{Le} for the values of the
$sl_2$ Witten-Reshetikhin-Turaev invariants \cite{Reshetikhin-Turaev}
\cite{Witten} of some particular $3$-manifolds, including the
Poincar\`e homology sphere, were expressed as infinite series similar
to~\eqref{eq:17}, and have well-defined values at the roots of unity.

The ring $\ZqN$ arises from the new invariant $I(M)$ of an
integral homology $3$-sphere $M$ that we announced
in~\cite{Habiro:rims2001} (see also~\cite{Ohtsuki:problemlist}).  (The
ring $\ZqN$ is denoted $\widehat{\modZ [q]}$ in~\cite{Habiro:rims2001}.)
The invariant $I(M)$ takes values in $\ZqN$ and unifies all the $sl_2$
Witten-Reshetikhin-Turaev
invariants $\tau _\zeta (M)$ defined at various roots $\zeta $ of unity; i.e.,
for any root $\zeta $ of unity we have
\[
  I(M)|_{q=\zeta }=\tau _\zeta (M).
\]
As we explained in \cite{Habiro:rims2001}, the existence of the invariant
$I(M)$ generalizes the previous integrality results \cite{Lawrence}
\cite{Murakami} \cite{Ohtsuki} \cite{Rozansky} on the
Witten-Reshetikhin-Turaev invariant for integral homology spheres.

The present paper was at first intended to provide the results on the
ring $\ZqN$ announced in~\cite{Habiro:rims2001} and those necessary
for the future papers \cite{Habiro:in-preparation} in which we will
prove the existence of $I(M)$.  From purely algebraic interests,
however, we will also study some generalizations of $\ZqN$ as follows.

Let $R$ be a commutative ring with unit, and let $R[q]$ denote the
polynomial ring over $R$ in an indeterminate $q$.  For each
$n\in \modN =\{1,2,\ldots\}$ let $\Phi _n(q)$ denote the $n$th cyclotomic
polynomial
\[
  \Phi _n(q)=\prod_{(i,n)=1}(q-\zeta ^i)\in \modZ [q],
\]
where $\zeta $ is a primitive $n$th root of unity.  If $S\subset \modN $ is a
subset, we set
\[
\PHI_S=\{\Phi _n(q)\zzzvert n\in S\}\subset \modZ [q],
\]
and let $\PHI_S^*$ denote the multiplicative set in $\modZ [q]$ generated
by $\PHI_S$, which we will regard as a directed set with respect to
the divisibility relation $|$.  The principal ideals $(f(q))\subset R[q]$
for $f(q)\in \PHI_S^*$ define a linear topology of the ring $R[q]$.
Define a commutative $R$-algebra $\RqS$ by
\begin{equation}
  \label{eq:1}
  \RqS = \varprojlim_{f(q)\in \PHI_S^*} R[q]/(f(q)),
\end{equation}
which we will call the {\em ($S$-)cyclotomic completion\/} of $R[q]$.
Since the sequence $(-1)^n(q)_n$, $n\ge 0$, is cofinal to the directed
set $\PHI_\modN ^*$, the definition~\eqref{eq:1} is consistent
with~\eqref{eq:3}.  Note that if $S$ is finite, then $\RqS$ is
identified with
the $(\prod\PHI_S)$-adic completion
of $R[q]$, where $\prod\PHI_S =
\prod_{f\in \PHI_S}f=\prod_{n\in S}\Phi _n(q)$.  In particular, we have
$\Rq{\{1\}}\simeq R[[q-1]]$ and $\Rq{\{2\}}\simeq R[[q+1]]$.  In
general, we have a natural isomorphism
\[
\ZqS\simeq\varprojlim_{S'\subset S, |S'|<S}\Zq{S'},
\]
where $S'$ runs through all the finite subsets of $S$.

We are interested in the behavior of natural homomorphisms among the
cyclotomic completions
$\RqS$ for various $R$ and $S$, and also those among the $\RqS$ and some
other rings.  First of all, if $g\zzzcolon R\rightarrow R'$ is a ring homomorphism, then
for each $S\subset \modN $ the homomorphism $g_q\zzzcolon R[q]\rightarrow R'[q]$, induced by $g$,
induces a ring homomorphism $g_S\zzzcolon \RqS\rightarrow \xq{R'}S$.  If $g$ is
injective (resp.  surjective), then so is $g_S$ (see
Lemma~\ref{lem:1}).

More interesting homomorphisms among cyclotomic completions
 are induced by inclusions
$S'\subset S\subset \modN $.  In this case, $\PHI_{S'}^*$ is a directed subset of
$\PHI_S^*$, and hence $\id_{R[q]}$ induces an $R$-algebra homomorphism
\[
\rho ^R_{S,S'}\zzzcolon \RqS\rightarrow \Rq{S'}.
\]
The rings $\RqS$ for $S\subset \modN $ and the homomorphisms $\rho ^R_{S,S'}$ form
a presheaf of rings over the set $\modN $ with the discrete topology;
i.e., we have $\Rq{\emptyset}=\{0\}$ and
$\rho ^R_{S,S''}=\rho ^R_{S',S''}\rho ^R_{S,S'}$ if $S''\subset S'\subset S\subset \modN $.

We will state a sufficient condition for $\rho ^R_{S,S'}$ to be injective
using a certain graph defined on the set $S$.  For each subset $S\subset \modN $, let
$\Gamma _R(S)$ denote the graph (with loop-edges) whose set of vertices
is $S$, and in which two elements $n,n'\in S$ are adjacent if and
only if either
\begin{enumerate}
\item $n=n'$,
\item $n/n'$ is an integer power of a prime $p$ such that $R$ is
$p$-adically separated, i.e., $\bigcap_{j\ge 0}p^jR=(0)$, or
\item $R=\{0\}$.
\end{enumerate}
If either one of the above conditions holds, then we write $n\Leftrightarrow _Rn'$.
For example, in $\Gamma _\modZ (\modN )$ two vertices $n,n'$ are adjacent if and
only if $n/n'$ is an integer power of a prime, and hence the graph
$\Gamma _\modZ (\modN )$ is connected; while the graph $\Gamma _\modQ (\modN )$ is discrete,
i.e., two distinct vertices are never adjacent.  Theorem~\ref{thm:17}
states that if $S'\subset S\subset \modN $ are subsets such that for any $n\in S$ there
is a sequence $S'\ni n'\Leftrightarrow _R\cdots\Leftrightarrow _Rn$ in $S$, then the homomorphism
$\rho ^R_{S,S'}$ is injective.
A nonempty subset
$S\subset \modN $ is said to be {\em $\Leftrightarrow _R$-connected} if the graph $\Gamma _R(S)$ is
connected.  It follows that
if $S\subset \modN $ is $\Leftrightarrow _R$-connected, then for any nonempty subset $S'\subset S$
the homomorphism $\rho ^R_{S,S'}$ is injective.
In particular, since $\modN $ is $\Leftrightarrow _\modZ $-connected, for each $n\in \modN $ the
homomorphism
\[
  \rho ^\modZ _{\modN ,\{n\}}\zzzcolon \ZqN\rightarrow  \Zq{\{n\}}
  \bigl(= \varprojlim_{j\ge 0} \modZ [q]/(\Phi _n(q)^j)\bigr)
\]
is injective.

If $\zeta $ is a primitive $n$th root of unity, then the
homomorphism
\[
\sigma ^\modZ _{\modN ,\zeta }\zzzcolon \ZqN\rightarrow \modZ [\zeta ][[q-\zeta ]],
\]
which is induced by the inclusion $\modZ [q]\subset \modZ [\zeta ][q]$ and factors
through $\rho ^\modZ _{\modN ,\{n\}}$, is injective (Theorem~\ref{thm:13}).  In
other words, each element of $\ZqN$ is uniquely determined by its
power series expansion in $q-\zeta $.  In particular, the invariant $I(M)$
of an integral homology sphere $M$ is completely determined by its
expansion in $q-\zeta $ for one root $\zeta $ of unity,
which in the case $\zeta =1$ is the Ohtsuki
series~\cite{Ohtsuki}.  Since $\modZ [\zeta ][[q-\zeta ]]$ is an integral domain,
it follows that so is $\ZqN$ (Corollary~\ref{thm:8}).

We are also interested in the homomorphism
\[
  \tau ^R_{S,T}\zzzcolon \RqS \rightarrow  P_T(R) = \prod_{n\in T}R[q]/(\Phi _n(q))
\]
for $T\subset S\subset \modN $, induced by the homomorphism $R[q]\rightarrow P_T(R)$,
$f(q)\mapsto(f(q)\mod(\Phi _n(q)))_{n\in T}$, where $R$ is a subring of the
field $\bar\modQ $ of algebraic numbers.  If $S$ is $\Leftrightarrow _R$-connected, and
for some $n\in S$ there are infinitely many elements $m\in T$ with
$m\Leftrightarrow _Rn$, then $\tau ^R_{S,T}$ is injective (Theorem~\ref{thm:3}).  In
particular, if $T\subset \modN $
contains infinitely many prime powers, then
$\tau ^\modZ _{\modN ,T}\zzzcolon \ZqN\rightarrow P_T(\modZ )$ is injective.  Hence it follows that if
$Z$ is a set of roots of unity containing infinitely many elements
of prime power order, then the homomorphism
\[
\tau ^R_{S,Z}\zzzcolon \ZqS \rightarrow  P_Z(\modZ ) = \prod_{\zeta \in Z}\modZ [\zeta ],
\]
induced by $\modZ [q]\rightarrow P_Z(\modZ )$, $f(q)\mapsto(f(\zeta ))_{\zeta \in Z}$, is
injective (Theorem~\ref{thm:2}).  In other words, each element in
$\ZqN$ is uniquely determined by its values at roots of unity in such
a set $Z$.  In particular, it follows that the invariant $I(M)$ of an
integral homology sphere $M$ is completely determined by the
Witten-Reshetikhin-Turaev invariants $\tau _\zeta (M)$ with $\zeta \in Z$.

Recall that a holomorphic function defined in a region is determined
either by the power series expansion at one point or by its values at
any infinitely many points contained in a compact set in the region.
The properties of $\ZqN$ described above may be interpreted that
$\ZqN$ behaves like the ring of ``holomorphic functions defined in
the set of the roots of unity''.  These properties are not as obvious as
they might first appear; the ring $\QqN$, which contains $\ZqN$ as a
subring, is quite contrasting.  We have an isomorphism
\[
\QqN \simeq \prod_{n\in \modN } \Qqn,
\]
see Section~\ref{sec:structure-qqh}.  It follows that
$\rho ^\modQ _{\modN ,\{n\}}$ for $n\in \modN $ and $\tau ^\modQ _{\modN ,\modN }$ are {\em not\/}
injective (but surjective), and that $\QqN$ is {\em not\/} an integral
domain.

The results stated above are more or less generalized in the later
sections.

The rest of the paper is organized as follows.  In
Section~\ref{sec:preliminaries} we fix some notations.
Section~\ref{sec:lemma-compl-polyn} deals with what we might call
``monic completions'' of $R[q]$, which are generalizations of
cyclotomic completions defined using monic polynomials instead of
cyclotomic polynomials.  In Section~\ref{sec:injectivity-rs-s} we
apply the results in Section~\ref{sec:lemma-compl-polyn} to cyclotomic
completions, and study the conditions for the homomorphisms
$\rho ^R_{S,S'}$ to be injective.  In Section~\ref{sec:z-cycl-compl} we
consider the power series expansion of the elements of $\RqS$ in
$q-\zeta $ with $\zeta \in R$ root of unity of order contained in $S$.  In
Section~\ref{sec:values-at-roots} we study the homomorphisms
$\tau ^R_{S,T}$ and $\tau ^R_{S,Z}$.  In Section~\ref{sec:some-remarks} we
give some remarks.

\section{Preliminaries}
\label{sec:preliminaries}

Throughout the paper, rings are unital and commutative, and
homomorphisms of rings are unital.  By ``homomorphism'' we will
usually mean a ring homomorphism.  Two rings that are considered to be
canonically isomorphic to each other will often be identified.  Also,
if a ring $R$ embeds into another ring $R'$ in a natural way, we will
often regard $R$ as a subring of $R'$.

If $R$ is a ring and $I\subset R$ is an ideal, then the $I$-adic completion
of $R$ will be denoted by
\[
R^I = \varprojlim_j R/I^j,
\]
and if $J\subset I$ is another ideal, then let
\[
\rho ^R_{J,I}\zzzcolon R^J\rightarrow R^I
\]
denote the homomorphism induced by $\id_R$.  These notation should not
cause confusions with $R[q]^S$ and $\rho ^R_{S,S'}$ defined in the
introduction.  We will further generalize these notations in the later
sections.  The ring $R$ is said to be {\em $I$-adically separated\/}
(resp. {\em $I$-adically complete\/}) if the natural homomorphism
$R\rightarrow R^I$ is injective (resp. an isomorphism).  Recall that $R$ is
$I$-adically separated if and only if $\bigcap_{j\ge 0}I^j=(0)$.

Let $\modN =\{1,2,\ldots\}$ denote the set of positive integers.  We
regard $\modN $ as a directed set with respect to the divisibility
relation $|$.  We will not use the letter $\modN $ for the same set
$\{1,2,\ldots\}$ when it is considered as an ordered set with the
usual order $\le $.

The letter $q$ will always denote an indeterminate.

\section{Monic completions of polynomial rings}
\label{sec:lemma-compl-polyn}

\subsection{Definitions and basic properties}
\label{sec:defin-basic-prop}

For a ring $R$ let $\modM _R$ denote the set of the monic polynomials in
$R[q]$, which is a directed set with respect to the divisibility
relation $|$.  For a subset $M\subset \modM _R$, let $M^*$ denote the
multiplicative set in $R[q]$ generated by $M$, which is a directed
subset of $\modM _R$.  The principal ideals $(f)$, $f\in M^*$, define a
linear topology of the ring $R[q]$, and let
\begin{equation}
  \label{eq:13}
  \RqM = \varprojlim_{f\in M^*} R[q]/(f)
\end{equation}
denote the completion.  (If $M=\{1\}$, then~\eqref{eq:13} implies
$\Rq{\{1\}}=R[q]/(1)=0$, which notationally contradicts to the
previous definition
$\Rq{\{1\}}=R[[q-1]]$.  In the rest of the paper, however,
``$\Rq{\{1\}}$'' will always mean $R[[q-1]]$.)

If $M'\subset M\subset \modM _R$, then $(M')^*$ is a directed subset of $M^*$, and
hence $\id_{R[q]}$ induces a homomorphism
\[
\rho ^R_{M,M'}\zzzcolon \RqM\rightarrow \Rq{M'}.
\]
We also extend the notation in the obvious way to
$\rho ^R_{M,I}\zzzcolon \RqM\rightarrow \Rq{I}$ for $M\subset \modM _R$ a subset and $I\subset R$ an ideal,
etc., if it is well defined.  (The general rule is that
$\rho ^R_{X,Y}\zzzcolon \Rq{X}\rightarrow \Rq{Y}$ is a homomorphism induced by $\id_{R[q]}$.)

If $M\subset \modM _R$ is finite, then the directed set $M^*$ is cofinal to the
sequence $(\prod M)^j$, $j\ge 0$.  Hence $\RqM$ is naturally isomorphic
to the $(\prod M)$-adic completion $\Rq{(\prod M)}$ of $R[q]$.  In
particular, if $f\in \modM _R$, then we have
\[
  \Rq{\{f\}}\simeq\Rq{(f)} = \varprojlim_{j}R[q]/(f)^j.
\]
If $M\subset \modM _R$ is infinite, then $\RqM$ is not an ideal-adic completion
in general, see for example Proposition~\ref{thm:9}.

If $M\subset \modM _R$, then the rings $\Rq{M'}$ for finite subsets $M'$ of $M$
and the natural homomorphisms $\rho ^R_{M',M''}$ for finite $M',M''$ with
$M''\subset M'\subset M$ form an inverse system of rings, of which the inverse
limit is naturally isomorphic to $\RqM$; i.e., we have
\begin{equation}
  \label{eq:7}
  \RqM \simeq \varprojlim_{M'\subset M,\,|M'|<\infty } \Rq{M'}.
\end{equation}

Let $h\zzzcolon R\rightarrow R'$ be a ring homomorphism.  Note that if $h$ is injective
(resp.  surjective), then so is the induced homomorphism
$h_q\zzzcolon R[q]\rightarrow R'[q]$.

\begin{lemma}
  \label{lem:1}
  Let $h\zzzcolon R\rightarrow R'$ be a ring homomorphism and let $M\subset \modM _R$ be at most
  countable.  If $h$ is injective (resp. surjective), then so is
  the homomorphism
  \[
  h_M\zzzcolon \RqM\rightarrow \xq{R'}{h(M)}
  \]
  induced by $h_q$.
\end{lemma}

\begin{proof}
  If $M\subset \{1\}$, then the result is trivial, so we assume this is not
  the case.  Since $M$ is at most countable, there is a (possibly
  terminating) sequence
  $g_0|g_1|\cdots$ in $M^*$, which is cofinal to $M^*$.  Note that the
  sequence $h(g_0)|h(g_1)|\cdots$, is cofinal to $h(M^*)=h(M)^*$.
  Since each $g_n$ is monic, each $a\in \RqM$ is uniquely expressed as
  an infinite sum
  \[
  a = \sum_{n\ge 0} a_n g_n,
  \]
  where $a_n\in R[q]$, $\deg a_n<\deg g_{n+1}-\deg g_n$ for $n\ge 0$.
  From this presentation of elements of $\RqM$ the result follows
  immediately.
\end{proof}

\subsection{Injectivity of the homomorphism $\rho ^R_{M,M'}$}
\label{sec:inject-homom}

Let $R$ be a ring, $I\subset R$ an ideal, and $f,g\in \modM _R$.  Let $\sqrt I$
denote the radical of $I$.  We write $f\overset{I}{\Rightarrow }_{R}g$, or
simply $f\overset{I}{\Rightarrow }g$, if $f\in \sqrt{(g)+I[q]}$, i.e., if
$f^m\in (g)+I[q]$ for some $m\ge 0$.  For $f,g\in \modM _R$, we write $f\Rightarrow _R g$,
or simply $f\Rightarrow g$, if we have $f\overset{I}{\Rightarrow }_R g$ for some ideal
$I\subset R$ with $\bigcap_{j\ge 0}I^j=(0)$.  Then $\Rightarrow _R$ defines a relation
on the set $\modM _R$.  Obviously, $g|f$ implies $f\Rightarrow g$.  Note also that
if $f\Rightarrow g$, $f|f'$, and $g'|g$, then $f'\Rightarrow g'$.

\begin{proposition}
  \label{thm:20}
  Let $R$ be a ring, and $f,g\in \modM _R$ with $f\Rightarrow _R g$.  Then the
  homomorphism $\rho ^R_{(fg),(f)}\zzzcolon \Rq{(fg)}\rightarrow \Rq{(f)}$ is injective.
\end{proposition}

\begin{proof}
  We first show that if $f\overset{I}{\Rightarrow }g$ and $R$ is $I$-adically
  complete, then $\rho ^R_{(fg),(f)}$ is an isomorphism.  Since $R\simeq
  R^I$ and $f$ is monic, we have
  \[
  \begin{split}
    &\Rq{(f)}
    \simeq\xq{R^I}{(f)}
    =\varprojlim_i(\varprojlim_j R/I^j)[q]/(f^i)\\
    &\quad \simeq\varprojlim_i(\varprojlim_j R[q]/((f^i)+I^j[q]))
    \simeq\Rq{(f)+I[q]}.
  \end{split}
  \]
  Similarly, $\Rq{(fg)}\simeq\Rq{(fg)+I[q]}$.  Since
  $f\overset{I}\Rightarrow g$, we have $((f)+I[q])^m\subset (f^m)+I[q]\subset (fg)+I[q]$ for
  some $m\ge 1$, while we obviously have $(fg)+I[q]\subset (f)+I[q]$.  Hence
  the $((f)+I[q])$-adic topology and the $((fg)+I[q])$-adic topology
  of $R[q]$ are the same.  Hence $\rho ^R_{(fg)+I[q],(f)+I[q]}$, which is
  identified with $\rho ^R_{(fg),(f)}$, is an isomorphism.

  Now consider the general case, where we have $f\overset{I}{\Rightarrow }_R g$
  and $R$ is $I$-adically separated.  We have a commutative diagram
  \[
  \begin{CD}
    \Rq{(fg)} @>{\rho ^R_{(fg),(f)}}>>\Rq{(f)} \\
    @VVV  @VVV \\
    \xq{R^I}{(fg)} @>>{\rho ^{R^I}_{(fg),(f)}}>\xq{R^I}{(f)}
  \end{CD}
  \]
  where vertical arrows are induced by the inclusion $R\subset R^I$, and
  hence are injective.
  Let $\bar I$ denote the closure of $I$ in $R^I$.
  Since $R^I$ is $\bar I$-adically complete and
  clearly $f\overset{\bar I}{\Rightarrow }_{R^I} g$, the above-proved case implies
  that $\rho ^{R^I}_{(fg),(f)}$ is an isomorphism.  Hence
  $\rho ^R_{(fg),(f)}$ is injective.
\end{proof}

For two subsets $M,M'\subset \modM _R$, we write $M'\prec M$ if $M'\subset M$ and for
each $f\in M$ there is a sequence $M'\ni f_0\Rightarrow f_1\Rightarrow \cdots \Rightarrow f_r=f$ in $M$.

Suppose that $M_0\prec M\subset \modM _R$.  Set
\[
\modF (M,M_0)=\{M'\subset M\zzzvert M_0\subset M',\, |M'\setminus M_0|<\infty  \},
\]
and
\[
\modF ^\prec(M,M_0)=\{M'\in \modF (M,M_0)\zzzvert M_0\prec M'\}\subset \modF (M,M_0).
\]
We will regard $\modF (M,M_0)$ as a directed set with respect to $\subset $, and
$\modF ^\prec(M,M_0)$ as a partially-ordered subset of $\modF (M,M_0)$.  Note that
if $M',M''\in \modF ^\prec(M,M_0)$ and $M''\subset M'$, then we have $M''\prec M'$.

\begin{lemma}
  \label{lem:4}
  If $M_0\prec M\subset \modM _R$, then $\modF ^\prec(M,M_0)$ is a cofinal directed
  subset of $\modF (M,M_0)$.
\end{lemma}

\begin{proof}
  It suffices to show that if $M'\in \modF (M,M_0)$, then there is
  $M''\in \modF ^\prec(M,M_0)$ with $M'\subset M''$.  For each $g\in M'\setminus M_0$ choose a
  sequence $M_0\ni g_0\Rightarrow \cdots\Rightarrow g_r=g$ in $M$ and set
  $U_g=\{g_1,\ldots,g_r\}$.  Set $M''=M_0\cup \bigcup_{g\in M'\setminus M_0}U_g$.  Then
  we have $M''\in \modF ^\prec(M,M_0)$ and $M'\subset M''$.
\end{proof}

\begin{theorem}
  \label{thm:21}
  If $R$ is a ring and $M_0\prec M\subset \modM _R$, then the homomorphism
  $\rho ^R_{M,M_0}\zzzcolon \RqM\rightarrow \Rq{M_0}$ is injective.
\end{theorem}

\begin{proof}
  By~\eqref{eq:7} and Lemma~\ref{lem:4} we have
  \[
  \RqM
  \simeq \varprojlim_{{M'}\in \modF (M,M_0)} \Rq{M'}
  \simeq \varprojlim_{{M'}\in \modF ^\prec(M,M_0)} \Rq{M'}.
  \]
  Hence it suffices to prove the theorem assuming that $M\setminus M_0$ is
  finite.  We can further assume that $|M\setminus M_0|=1$.  Let $g\in M\setminus M_0$ be
  the unique element.

  First we assume that $M_0=\{f_1,\ldots,f_n\}$ ($n\ge 1$) is
  finite.  Set $f=f_1\cdots f_n$.  Since $f_i\Rightarrow g$ for some $i\in \{1,\ldots,n\}$,
  we have $f\Rightarrow g$.
  By Proposition~\ref{thm:20},
  $\rho ^R_{(fg),(f)}$ is injective.
  Since $\Rq{M_0}=\Rq{(f)}$ and
  $\RqM=\Rq{(fg)}$, it follows that $\rho ^R_{M,M_0}$ is injective.

  Now assume that $M_0$ is infinite.  Choose an element
  $g_0\in M_0$ with $g_0\Rightarrow g$.  We have $\Rq{M_0} \simeq
  \varprojlim_{U\in \modF (M_0,\{g_0\})}\Rq{U}$ and $\RqM \simeq
  \varprojlim_{U\in \modF (M_0,\{g_0\})}\Rq{U\cup \{g\}}$.  For each
  ${U\in \modF (M_0,\{g_0\})}$ we have $U\prec U\cup \{g\}$.  Hence it follows
  from the above-proved case that the homomorphism
  $\rho ^R_{U\cup \{g\},U}\zzzcolon \Rq{U\cup \{g\}}\rightarrow \Rq U$ is injective.
  Since $\rho ^R_{M,M_0}$ is the
  inverse limit of the $\rho ^R_{U\cup \{g\},U}$ for $U\in \modF (M_0,\{g_0\})$,
  it is injective.
\end{proof}

A subset $M\subset \modM _R$ is said to be {\em $\Rightarrow _R$-connected\/} if $M$ is
not empty and for each $f,f'\in M$ there is a sequence $f=f_0\Rightarrow _R
f_1\Rightarrow _R\cdots\Rightarrow _R f_r=f'$ ($r\ge 0$) in $M$.  Note that if $M$ is
$\Rightarrow _R$-connected, then for any nonempty subset $M'\subset M$ we have
$M'\prec M$.  The following follows immediately from
Theorem~\ref{thm:21}.

\begin{corollary}
  \label{thm:14}
  If $R$ is a ring, and $M\subset \modM _R$ is a $\Rightarrow _R$-connected subset, then
  for any nonempty subset $M'\subset M$ the homomorphism
  $\rho ^R_{M,M'}\zzzcolon \RqM\rightarrow \Rq{M'}$ is injective.
\end{corollary}

\section{Injectivity of $\rho ^R_{S,S'}$}
\label{sec:injectivity-rs-s}

If $R$ a ring, and $S\subset \modN $ is a subset, then we have
$\RqS=\Rq{\PHI_S}$.  If $S'\subset S$, then we have
\[
\rho ^R_{S,S'}=\rho ^R_{\PHI_S,\PHI_{S'}}\zzzcolon \RqS\rightarrow \Rq{S'}.
\]

We will use the following well-known properties of cyclotomic
polynomials.
\begin{lemma}
  \label{lem:3}
  (1) Let $n\in \modN $, $p$ a prime, and $e\ge 1$.
  Then we have
  \begin{equation}
    \label{eq:6}
    \Phi _{p^en}(q)\equiv \Phi _n(q)^d \pmod{(p)},
  \end{equation}
  in $\modZ [q]$, where $d=\deg\Phi _{p^en}(q)/\deg\Phi _n(q)$.  (We have
  $d=(p-1)p^{e-1}$ if $(n,p)=1$ and $d=p^e$ if $p|n$.)

  (2) If $m,n\in \modN $, and $n/m\in \modQ $ is not an integer power of a
  prime, then we have $(\Phi _n(q),\Phi _m(q))=(1)$ in $\modZ [q]$.
\end{lemma}

For $m,n\in \modN $, we define $c_{m,n}\in \{0,1\}\cup \{p\zzzvert \text{prime}\}$ by
\begin{enumerate}
\item $c_{n,n}=0$,
\item $c_{m,n}=p$ if $p$ is a prime and $n/m=p^j$ for some
  $j\in \modZ \setminus \{0\}$, and
\item $c_{m,n}=1$ if $n/m$ is not an integer power of a prime.
\end{enumerate}
Note that $c_{m,n}=c_{n,m}$ for all $m,n\in \modN $.  It is straightforward
to see that $m\Leftrightarrow _Rn$ if and only if $R$ is $(c_{m,n})$-adically
separated.

Lemma~\ref{lem:3} implies that for each $m,n\in \modN $ we have
$\Phi _m(q)\in \sqrt {(\Phi _n(q),c_{m,n})}$ in $R[q]$, i.e.,
$\Phi _m(q)\overset{(c_{m,n})}\Rightarrow _R\Phi _n(q)$.  It follows that if $m\Leftrightarrow _Rn$,
then we have $\Phi _m(q)\Rightarrow _R\Phi _n(q)$.
Note also that if $S\subset \modN $ is
$\Leftrightarrow _R$-connected, then $\PHI_S$ is $\Rightarrow _R$-connected.  The following
follows immediately from Theorem~\ref{thm:21} and
Corollary~\ref{thm:14}.

\begin{theorem}
  \label{thm:17}
  Let $R$ be a ring and let $S'\subset S\subset \modN $.  Suppose that each connected
  component of the graph $\Gamma _R(S)$ contains at least one vertex of
  $\Gamma _R(S')$.  (In other words, for each element $n\in S$, there is a
  sequence $S'\ni n'\Leftrightarrow _R\cdots\Leftrightarrow _Rn$ in $S$.)  Then the homomorphism
  $\rho ^R_{S,S'}$ is injective.

  In particular, if $S\subset \modN $ is $\Leftrightarrow _R$-connected, then for any nonempty
  subset $S'\subset S$ the homomorphism $\rho ^R_{S,S'}\zzzcolon \RqS\rightarrow \Rq{S'}$ is
  injective.  More particularly, for any nonempty subset $S'\subset \modN $ the
  homomorphism
  $\rho ^\modZ _{\modN ,S'}\zzzcolon \ZqN\rightarrow \Zq{S'}$ is injective.
\end{theorem}

We remark that the special case of Theorem~\ref{thm:17} where $R=\modZ $,
$S=\modN $, and $S'=\{1\}$ is obtained also by P.~Vogel.  Another proof of
a special case of Theorem~\ref{thm:17} is sketched in
Remark~\ref{thm:11}.

For each $n\in \modN $ set $\langle n\rangle =\{m\in \modN \zzzvert m|n\}$.  Since
$\prod\PHI_{\langle n\rangle }=\prod_{m|n}\Phi _m(q)=q^n-1$, we have
\[
\Rq{\langle n\rangle } =\Rq{(q^n-1)} = \varprojlim_j R[q]/(q^n-1)^j.
\]
Note that the set $\langle n\rangle $ is $\Leftrightarrow _R$-connected if and only if for each
prime factor $p$ of $n$ the ring $R$ is $p$-adically separated.
A subset $S\subset \modN $ will be called {\em $R$-admissible} if $S$ is a
$\Leftrightarrow _R$-connected, directed subset of $\modN $ such that $n\in S$ implies
$\langle n\rangle \subset S$.  Note that a subset $S\subset \modN $ is finite and $R$-admissible
if and only if there is $n\in \modN $ such that $S=\langle n\rangle $ and $R$ is
$p$-adically separated for each
prime factor $p$ of $n$.  Note also
that an $R$-admissible subset $S\subset \modN $ satisfies
$S=\bigcup_{n\in S}\langle n\rangle $, and hence we have
$\RqS\simeq\varprojlim_{n\in S}\Rq{\langle n\rangle }$.
The following easily follows from Theorem~\ref{thm:17}.

\begin{corollary}
  \label{thm:16}
  Let $R$ be a ring, and let $S\subset \modN $ be $R$-admissible.
  Then for each $m,n\in S$ with
  $m|n$ the homomorphism $\rho ^R_{\langle n\rangle ,\langle m\rangle }\zzzcolon \Rq{\langle n\rangle }\rightarrow \Rq{\langle m\rangle }$ is
  injective.  Hence $\RqS$ can be regarded as the intersection
  $\bigcap_{n\in S}\Rq{\langle n\rangle }$, where the $\Rq{\langle n\rangle }$, $n\in S$, are
  regarded as $R$-subalgebras of $\Rq{\langle 1\rangle }=R[[q-1]]$.

  In particular, if $m,n\in \modN $ and $m|n$, then
  $\rho ^\modZ _{\langle n\rangle ,\langle m\rangle }\zzzcolon \Zq{\langle n\rangle }\rightarrow \Zq{\langle m\rangle }$ is injective.  We have
  $\ZqN=\bigcap_{n\in \modN }\Zq{\langle n\rangle }$.
\end{corollary}

We will see in Proposition~\ref{thm:10} that if $m|n$ and $m\neq n$,
then $\rho ^\modZ _{\langle n\rangle ,\langle m\rangle }$ is not surjective.

\section{Expansions at roots of unity}
\label{sec:z-cycl-compl}

For an integral domain $R$ of characteristic $0$ let $Z^R$ denote the
set of the roots of unity in $R$.  If $S\subset \modN $, then set
$Z^R_S=\{\zeta \in Z^R\zzzvert \ord\zeta \in S\}$.  For a subset $Z\subset Z^R$ set
\[
  \RqZ = \Rq{M_Z},
\]
where $M_Z = \{q-\zeta \zzzvert \zeta \in Z\}\subset \modM _R$.  If $Z'\subset Z$, then set
\[
\rho ^R_{Z,Z'}=\rho ^R_{M_Z,M_{Z'}}\zzzcolon \RqZ\rightarrow \Rq{Z'}.
\]
(Although we have $1\in Z$ and $1\in \modN $, the notation $\Rq{\{1\}}$ is not
ambiguous because $1$ is the unique primitive $1$st root of unity.)

For a subset $Z\subset Z^R$ set $N_Z=\{\ord\zeta \zzzvert \zeta \in Z\}$, and in particular
set $N_R=N_{Z^R}$.  If $S\subset N_R$, then we have
\[
  \RqS \simeq \Rq{Z^R_S}.
\]

\begin{lemma}
  \label{lem:2}
  Let $R$ be an integral domain of characteristic $0$, and let
  $\zeta ,\zeta '\in Z^R$.  Then the following conditions are equivalent.
  \begin{enumerate}
  \item $(q-\zeta )\Rightarrow _R(q-\zeta ')$,
  \item $R$ is $(\zeta -\zeta ')$-adically separated,
  \item $\ord(\zeta ^{-1}\zeta ')$ is a power of some prime $p$ such that $R$ is
    $p$-adically separated.
  \end{enumerate}
\end{lemma}

\begin{proof}
  If (1) holds, then we have $(q-\zeta )^m\in (q-\zeta ')+I[q]$ for some $m\ge 0$
  and $R$ is $I$-adically separated.  It follows that $(\zeta '-\zeta )^m\in I$,
  and hence $R$ is $(\zeta '-\zeta )$-adically separated.  Hence we have (2).

  The other implications (2)$\Rightarrow $(1) and (2)$\Leftrightarrow $(3) are
  straightforward.
\end{proof}

Let $\Leftrightarrow _R$ denote the relation on $Z^R$ such that for $\zeta ,\zeta '\in Z^R$ we
have
$\zeta \Leftrightarrow _R\zeta '$ if and only if either one of the conditions in
Lemma~\ref{lem:2} holds.  The following follows immediately from
Corollary~\ref{thm:14}.

\begin{theorem}
  \label{thm:1}
  Let $R$ be an integral domain of characteristic $0$ and let $Z\subset Z^R$
  be a $\Leftrightarrow _R$-connected subset.  Then for any nonempty subset $Z'\subset Z$
  the homomorphism $\rho ^R_{Z,Z'}\zzzcolon \RqZ\rightarrow \Rq{Z'}$ is injective.
\end{theorem}

If $\zeta ,\zeta '\in Z^R$, then $\zeta \Leftrightarrow _R\zeta '$ implies $\ord\zeta \Leftrightarrow _R\ord\zeta '$.  (The
converse, however, does not holds.)  It follows that if $Z\subset Z^R$ is
$\Leftrightarrow _R$-connected, then $N_Z$ is $\Leftrightarrow _R$-connected.

\begin{remark}
  \label{thm:11}
  We sketch below another proof using Theorem~\ref{thm:1} of
  the special case of Theorem~\ref{thm:17} where $S$ is
  $\Leftrightarrow _R$-connected and $R$ is an
  integral domain of characteristic $0$ such that $R$ is $p$-adically
  separated for any prime $p$.  Let $k$ be the quotient field of $R$
  and let $\bar k$ be the algebraic closure of $k$.  Let $\tilde
  R\subset \bar{k}$ be the $R$-subalgebra generated by the elements of
  $Z_S^{\bar k}$.  In view of Lemma~\ref{lem:1}, it suffices to see
  that $\rho ^{\tilde R}_{S,S'}$ is injective.  Since $\Leftrightarrow _R$-connectivity
  of $S$ implies that of $Z_S$, the homomorphism $\rho ^{\tilde
    R}_{S,S'}=\rho ^{\tilde R}_{Z_S,Z_{S'}}$ is injective by
  Theorem~\ref{thm:1}.
\end{remark}

\begin{theorem}
  \label{thm:13}
  Let $R$ be an integral domain of characteristic $0$, $S\subset \modN $ a
  $\Leftrightarrow _R$-connected subset, and $n\in S$.  Assume that $R$ is
  $p$-adically separated for each odd prime factor $p$ of $n$, and
  also that if $4|n$, then $R$ is $2$-adically
  separated.  Let $\zeta $ be a primitive $n$th root of unity, which may
  or may not be contained in $R$.  Then the homomorphism
  \[
  \sigma ^R_{S,\zeta }\zzzcolon \RqS\rightarrow R[\zeta ][[q-\zeta ]]
  \]
  induced by $R[q]\subset R[\zeta ][q]$ is injective.

  In particular, for any root $\zeta $ of unity the homomorphism
  $\sigma ^\modZ _{\modN ,\zeta }\zzzcolon \ZqN\rightarrow \modZ [\zeta ][[q-\zeta ]]$ is injective.
\end{theorem}

\begin{proof}
  The homomorphism $\sigma ^R_{S,\zeta }$ is the composition of the following
  three homomorphisms
  \[
  \begin{CD}
  \RqS
  @>\rho ^R_{S,\{n\}}>> \Rqn
  @>i>> \xq{R[\zeta ]}{\{n\}}
  @>\rho ^{R[\zeta ]}_{\{n\},(q-\zeta )}>> R[\zeta ][[q-\zeta ]],
  \end{CD}
  \]
  the first two arrows of which are injective by Theorem~\ref{thm:17}
  and Lemma~\ref{lem:1}, respectively.  Hence it suffices to prove
  that $\rho ^{R[\zeta ]}_{\{n\},(q-\zeta )}$ is injective.  We may assume
  $\zeta \in R$, hence $R=R[\zeta ]$.

  For each $m$ with $m|n$, set $Z_m=Z^R_{\{m\}}=\{\zeta \in Z^R\zzzvert \ord\zeta =m\}$.
  By $\Rqn\simeq\Rq{Z_n}$ and
  Theorem~\ref{thm:1}, it suffices to prove that the set $Z_n$ is
  $\Leftrightarrow _R$-connected.  The case $n=1$ is trivial, so we assume not.
  Let $n=p_1^{e_1}\cdots p_r^{e_r}$ be a factorization into
  prime powers, where $p_1,\ldots,p_r$ are distinct primes and
  $e_1,\ldots,e_r\ge 1$.  There is a bijection
  \[
  Z_{p_1^{e_1}}\times \cdots \times Z_{p_r^{e_r}}\simeqto Z_n,\quad
  (\xi _1,\ldots,\xi _r)\mapsto \xi _1\cdots \xi _r.
  \]
  It suffices to show that if
  $(\xi _1,\ldots,\xi _r),(\xi '_1,\ldots,\xi '_r)\in Z_{p_1^{e_1}}\times \cdots\times Z_{p_r^{e_r}}$
  satisfies $\xi _j=\xi '_j$ for all $j\in \{1,\ldots,r\}\setminus \{i\}$ and
  $\xi _i\neq\xi '_i$ for some $i$, then we have $\xi _1\cdots\xi _r\Leftrightarrow _R
  \xi '_1\cdots \xi '_r$, which is equivalent to that $\xi _i\Leftrightarrow _R\xi '_i$.
  Since $Z_2=\{-1\}$ contains only one element, the case $p_i=2$ and
  $e_i=1$
  does not occur.  We have
  $(\xi _i-\xi '_i)\subset \sqrt{(p_i)}$, and hence $\xi _i\Leftrightarrow _R\xi '_i$.
\end{proof}

\begin{corollary}
  \label{thm:8}
  Let $R$ be an integral domain of characteristic $0$, and $S\subset \modN $ a
  $\Leftrightarrow _R$-connected subset.  Suppose that there is $n\in S$ such that
  $R$ is $p$-adically separated for each odd prime factor $p$ of
  $n$, and if $4|n$, then $R$ is also $2$-adically separated.  Then
  the ring $\RqS$ is an integral domain.

  In particular, $\ZqS$ is an integral domain for any nonempty subset
  $S\subset \modN $.
\end{corollary}

\begin{proof}
  The result follows from Theorem~\ref{thm:13} and the fact that the
  formal power series ring $R[\zeta ][[q-\zeta ]]$ is an integral domain.
\end{proof}

\section{Values at roots of unity}
\label{sec:values-at-roots}

\begin{theorem}
  \label{thm:3}
  Let $R$ be a subring of the field $\bar\modQ $ of algebraic numbers,
  $S\subset \modN $ a $\Leftrightarrow _R$-connected subset, and $T\subset S$ a subset.  Suppose
  that for some $n\in S$ there are infinitely many elements $m\in T$ with
  $m\Leftrightarrow _Rn$.  Then the homomorphism $\tau ^R_{S,T}\zzzcolon \RqS\rightarrow P_T(R)$ is
  injective.

  In particular, if $R$ is a subring of the ring of algebraic
  integers, then, for any subset $T\subset \modN $ containing infinitely many
  prime powers, $\tau ^R_{\modN ,T}\zzzcolon \RqN\rightarrow P_T(R)$ is
  injective.
\end{theorem}

\begin{proof}
  Suppose for contradiction that there is a nonzero element $a\in \RqS$
  with $\tau ^R_{S,T}(a)=0$.  By Theorem~\ref{thm:17}, $\rho ^R_{S,\{n\}} $
  is injective, and therefore we have $\rho ^R_{S,\{n\}}(a)\neq 0$.
  Hence we can write $\rho ^R_{S,\{n\}}(a)=\sum_{j=l}^\infty a_j\Phi _n(q)^j$,
  where $l\ge 0$ and $a_j\in R[q]$ for $j\ge l$ with $a_l\not\in (\Phi _n(q))$.
  There are infinitely many elements $m_1,m_2,\ldots\in T$ with
  $m_i\Leftrightarrow _Rn$ and $n|m_i$.  For each $i$, $m_i/n$ is a power of a prime
  $p_i$ such that $R$ is $p_i$-adically separated.  It follows from
  $\tau ^R_{S,T}(a)=0$ that $\Phi _{m_i}(q)|a$ in $\RqS$ for each $i$.
  Since $R$ is an integral domain of characteristic $0$, we can show
  by induction that $\Phi _{m_1}(q)\cdots\Phi _{m_k}(q)|a$ in $\RqS$ for
  each $k\ge 0$, and hence we have $\Phi _{m_1}(q)\cdots
  \Phi _{m_k}(q)|\rho ^R_{S,\{n\}}(a)$ in $\Rqn$.  By~\eqref{eq:6} we have
  $\Phi _{m_i}(q)\in (p_i,\Phi _n(q))$ for each $i$.  Hence we have
  $\Phi _{m_1}(q)\cdots \Phi _{m_k}(q)\in (p_1\cdots p_k,\Phi _n(q))$.  In other
  words, for each $k\ge 0$, $\bar a_l= a_l\mod{(\Phi _n(q))}
  \in R[q]/(\Phi _n(q))$ is divisible by $p_1\cdots p_k$.  Note that
  $R[q]/(\Phi _n(q))=R\oplus Rq\oplus\cdots\oplus Rq^{d-1}$ with
  $d=\deg\Phi _n(q)$, and $\bar a_l$ is expressed as a polynomial in $q$
  of degree$<d$, each coefficient of which is divisible by $p_1\cdots
  p_k$ in $R$ for $k\ge 0$.  Since $R$ is a subring of $\bar \modQ $ and
  each $p_i$ is a non-unit in $R$, it follows that the coefficients of
  $\bar a_l$ are zero.  Consequently, we have $a_l\in (\Phi _n(q))$.
\end{proof}

\begin{proposition}
  \label{thm:9}
  Let $R$ be a subring of $\bar\modQ $, and $S\subset \modN $ an infinite subset.
  Then the completion $\RqS$ of $R[q]$ is {\em not} an ideal-adic
  completion, i.e., there is no ideal $I$ in $R[q]$ such that
  $\id_{R[q]}$ induces an isomorphism
  $\RqS\simeq\varprojlim_jR[q]/I^j$.
\end{proposition}

\begin{proof}
  Let $I\subset R[q]$ is a nonzero ideal, and $f(q)\in I$ a nonzero element.
  Then there are only finitely many elements $n\in S$ with $\Phi _n(q)|f(q)$.
  For each $n\in R$ the power $f(q)^j$ for $j\ge 1$ is divisible by
  $\Phi _n(q)$ if and only if $f(q)$ is divisible by $\Phi _n(q)$.  It
  follows that $f(q)^j$ does not converge to $0$ as $j\rightarrow 0$ in $R[q]$
  with the topology defining the completion $\RqS$.  Hence we have
  $\RqS\not\simeq\varprojlim_{j}R[q]/I^j$.
\end{proof}

Let $R$ be a subring of $\bar\modQ $, and let $Z\subset Z^{\bar \modQ }$ be a
subset.  Set
\[
  P_Z(R) = \prod_{\zeta \in Z}R[\zeta ],
\]
which generalize the definition of $P_Z(\modZ )$.
If $S\subset \modN $ is a subset and $Z\subset Z_S^{\bar\modQ }$, then let
\[
  \tau ^R_{S,Z}\zzzcolon \RqS\rightarrow P_Z(R)
\]
denote the homomorphism induced by $R[q]\rightarrow P_Z(R)$,
$f(q)\mapsto(f(\zeta ))_{\zeta \in Z}$.

\begin{theorem}
  \label{thm:2}
  Let $R$ be a subring of $\bar\modQ $, and let $S\subset \modN $ and
  $Z\subset Z^{\bar\modQ }_S$ be subsets.  Suppose that there is an element
  $n\in S$ such that for infinitely many $\zeta \in Z$ we have $\ord\zeta \Leftrightarrow _R n$.
  Then the homomorphism $\tau ^R_{S,Z}\zzzcolon \RqS\rightarrow P_Z(R)$ is injective.

  In particular, if $R$ is a subring of the ring of algebraic
  integers, and $Z\subset Z^{\bar\modQ }$ is a subset containing infinitely many
  elements of prime power orders, then $\tau ^R_{S,Z}\zzzcolon \RqS\rightarrow P_Z(R)$ is
  injective.
\end{theorem}

\begin{proof}
  Set $N_Z=\{\ord\zeta \zzzvert \zeta \in Z\}\subset \modN $.  Let $\gamma \zzzcolon P_{N_Z}(R)\rightarrow P_Z(R)$ be
  the homomorphism defined by
  $\gamma ((f_n(q))_{n\in N_Z})=(f_{n_\zeta }(\zeta ))_{\zeta \in Z}$.  Since $\gamma $ is the
  direct product of the injective homomorphisms
  $R[q]/(\Phi _n(q))\rightarrow \prod_{\zeta \in Z, \ord\zeta =n}R[\zeta ]$,
  $f(q)\mapsto(f(\zeta ))_\zeta $, it follows that $\gamma $ is injective.  We have
  $\tau ^R_{S,Z}=\gamma \tau ^R_{S,N_Z}$, where $\tau ^R_{S,N_Z}\zzzcolon \RqS\rightarrow P_{N_Z}(R)$
  is injective by Theorem~\ref{thm:3}.  Hence $\tau ^R_{S,Z}$ is
  injective.
\end{proof}

\begin{conjecture}
  \label{thm:5}
  For any infinite subset $Z\subset Z^{\bar\modQ }$, the homomorphism
  $\tau ^\modZ _{\modN ,Z}\zzzcolon \ZqN\rightarrow P_Z(\modZ )$ is injective.
\end{conjecture}

If $Z'\subset Z\subset Z^R$, then we have a homomorphism
\[
  \tau ^R_{Z,Z'}\zzzcolon \RqZ\rightarrow  P_Z(R),
\]
induced by $R[q]\rightarrow P_Z(R)$, $f(q)\mapsto (f(\zeta ))_\zeta $.

\begin{theorem}
  \label{thm:15}
  Let $R$ be a subring of $\bar \modQ $, let $Z\subset Z^R$ a $\Leftrightarrow _R$-connected
  subset, and let $Z'\subset Z$.  Suppose that for some $\zeta \in Z$ there are
  infinitely many elements $\xi \in Z'$ with $\xi \Leftrightarrow _R\zeta $.  Then the
  homomorphism $\tau ^R_{Z,Z'}\zzzcolon \RqZ\rightarrow P_{Z'}(R)$ is injective.
\end{theorem}

\begin{proof}
  The proof is similar to that of Theorem~\ref{thm:3} with the
  cyclotomic polynomials replaced with the polynomials $q-\zeta $ with
  $\zeta $ a roots of unity.  The details are left to the reader.
\end{proof}

\section{Remarks}
\label{sec:some-remarks}

\subsection{Units in $\ZqS$}
\label{sec:invertibility-q}

If $R$ is a ring and $S\subset \modM _R$ is a subset consisting of monic
polynomials with the
constant terms being units in $R$, then
the element $q$ is invertible in $\RqS$.  In particular, we have an
explicit formula for $q^{-1}\in \RqN$ as follows.

\begin{proposition}
  \label{thm:6}
  For any ring $R$ the element $q\in \RqN$ is invertible with the
  inverse
  \[
  q^{-1} = \sum_{n\ge 0}q^n(q)_n.
  \]
\end{proposition}

\begin{proof}
  \(
  q\sum_{n\ge 0}q^n(q)_n
  =\sum_{n\ge 0}q^{n+1}(q)_n
  =\sum_{n\ge 0}(1-(1-q^{n+1}))(q)_n
  =\sum_{n\ge 0}((q)_n-(q)_{n+1})
  =(q)_0
  =1.
  \)
\end{proof}

For each subset $S\subset \modN $ the inclusion $\modZ [q]\subset \modZ [q,q^{-1}]$ induces an
isomorphism
\[
\ZqS \simeq \varprojlim_{f\in \PHI_S^*}\modZ [q,q^{-1}]/(f),
\]
via which we will identify these two rings.  If $S\neq\emptyset$, then,
since $\bigcap_{f\in \PHI_S^*}(f)=(0)$ in $\modZ [q,q^{-1}]$, the natural
homomorphism $\modZ [q,q^{-1}]\rightarrow \ZqS$ is injective and regarded as
inclusion.

For a ring $R$ let $U(R)$ denote the (multiplicative) group of the units
in $R$.  If $S\neq\emptyset$, then we have
\[
U(\modZ [q,q^{-1}])\subset U(\ZqN).
\]
It is well known that $U(\modZ [q,q^{-1}])=\{\pm q^i\zzzvert i\in \modZ \}$.
If we regard $\ZqN$ and the $\Zq{\langle n\rangle }$ as subrings of
$\Zq{\langle 1\rangle }=\modZ [[q-1]]$ as in Corollary~\ref{thm:16}, then we have
\[
U(\ZqN)=\bigcap_{n\in \modN }U(\Zq{\langle n\rangle }).
\]

\begin{conjecture}
  \label{thm:7}
  We have $U(\ZqN)=\{\pm q^i\zzzvert i\in \modZ \}$.
\end{conjecture}

\begin{remark}
  One might expect that Conjecture~\ref{thm:7} would
  generalize to any infinite, $\modZ $-admissible subset $S\subset \modN $,
  but this is not the case.  For odd $m\ge 3$
  consider the element $\gamma _m=\sum_{i=0}^{m-1}(-1)^iq^i\in \modZ [q]$, which
  is known to define a unit in the ring $\modZ [q]/(q^n-1)$ with
  $(n,2m)=1$ and is called an ``alternating unit'', see
  \cite{Sehgal}.  For such $n$, it follows that
  there are $u,v\in \modZ [q]$ such that
  $\gamma _mu=1+v\Phi _n(q)$.  Since $1+v\Phi _n(q)$ is a unit in
  $\Zq{\langle n\rangle }$, it follows that $\gamma _m$ is a unit in $\Zq{\langle n\rangle }$.
  Set $S=\{n\in \modN \zzzvert (n,2m)=1\}$.  Then it is straightforward to check
  that $\gamma _m$
  defines a unit in $\ZqS$ (hence also in $\Zq{S'}$ for any $S'\subset S$).
  Consequently, we have $U(\ZqS)\subsetneq\{\pm q^i\zzzvert i\in \modZ \}$.
\end{remark}

\subsection{A localization of $\ZqN$}
\label{sec:some-localizations}

In some application~\cite{Habiro:in-preparation}, it will be natural
to consider the following type of localization of $\ZqN$.  Recall from
Proposition~\ref{thm:8} that $\ZqN$ is an integral domain.  Let
$Q(\ZqN)$ denote the quotient field of $\ZqN$.  We will consider the
the $\ZqN$-subalgebra $\ZqN[\PHI_\modN ^{-1}]$ of $Q(\ZqN)$ generated by the
elements $\Phi _n(q)^{-1}$ for $n\in \modN $.  Alternatively,
$\ZqN[\PHI_\modN ^{-1}]$ may be defined as the subring of $Q(\ZqN)$
consisting of the  fractions $f(q)/g(q)$ with $f(q)\in \ZqN$ and
$g(q)\in \PHI_\modN ^*$.  
Similarly, let $\modZ [q,q^{-1}][\PHI_\modN ^{-1}]$ denote the
$\modZ [q,q^{-1}]$-subalgebra of the quotient field $\modQ (q)(\subset Q(\ZqN))$ of
$\modZ [q,q^{-1}]$ generated by the elements  $\Phi _n(q)^{-1}$ for $n\in \modN $,
which may alternatively defined as the subring of $\modQ (q)$ consisting
of the  fractions $f(q)/g(q)$ with $f(q)\in \modZ [q,q^{-1}]$ and
$g(q)\in \PHI_\modN ^*$.  

\begin{proposition}
  \label{thm:18}
  We have $\ZqN[\PHI_\modN ^{-1}]=\ZqN + \modZ [q,q^{-1}][\PHI_\modN ^{-1}]$.
\end{proposition}

\begin{proof}
  The inclusion $\supset $ is obvious; we will show the other inclusion.
  Since
  \[
  \ZqN[\PHI_\modN ^{-1}]
  =\bigcup_{f(q)\in \PHI_\modN ^*}\frac1{f(q)}\ZqN,
  \]
  it suffices to show that for each $f(q)\in \PHI_\modN ^*$ we have
  \[
  \frac1{f(q)}\ZqN \subset \ZqN + \frac1{f(q)}\modZ [q,q^{-1}].
  \]
  By multiplying $f(q)$, we need to show that
  \[
  \ZqN \subset f(q)\ZqN + \modZ [q,q^{-1}],
  \]
  which follows from $\ZqN\simeq
  \varprojlim_{g(q)\in \PHI_\modN ^*}\modZ [q,q^{-1}]/(f(q)g(q))$. 
\end{proof}

\begin{proposition}
  \label{thm:19}
  We have
  \[
  \ZqN\cap \modZ [q,q^{-1}][\PHI_\modN ^{-1}] = \modZ [q,q^{-1}].
  \]
\end{proposition}

\begin{proof}
  The inclusion $\supset $ is obvious; we will show the other inclusion.
  Suppose that $f(q)=g(q)/h(q)\in \ZqN\cap \modZ [q,q^{-1}][\PHI_\modN ^{-1}]$,
  where $g(q)\in \modZ [q,q^{-1}]$ and $h(q)\in \PHI_\modN ^*$, and $g(q)$ and
  $h(q)$ are coprime.  We will show that $f(q)\in \modZ [q,q^{-1}]$, i.e.,
  $h(q)|g(q)$ in $\modZ [q,q^{-1}]$.  Assume for contradiction that
  $h(q)\neq1$.  Choose $n\in \modN $ such that $\Phi _n(q)|h(q)$.  We have
  $g(q)=f(q)h(q)$.  Since $\Phi _n(q)|h(q)$, we have $\Phi _n(q)|g(q)$,
  which is a contradiction.  Hence we have $h(q)=1$, and we obviously
  have $h(q)|g(q)$.
\end{proof}

\subsection{Modules}
\label{sec:modules}

We can define cyclotomic completions also for any $\modZ $-modules as
follows.  Let $A$ be a $\modZ $-module, and let $A[q]$ denote the
$\modZ [q]$-module of polynomials in $q$ with coefficients in $A$.  For
each $S\subset \modN $ let $\AqS$ denote the completion
\[
  \AqS= \varprojlim_{f\in \PHI_S^*}A[q]/ f A[q].
\]
If $A$ is a ring, then this definition of $\AqS$ is compatible with
the previous one.  Some results in the present paper can be generalized to the
$\AqS$.

For example, Theorem~\ref{thm:17} is generalized as follows.
Let $\Leftrightarrow _A$ denote the relation on $\modN $
such that $m\Leftrightarrow _An$ if and only if either we have $A=0$, or $m/n$ is
an integer power of a prime $p$ with $A$ being $p$-adically separated.

\begin{theorem}
  \label{thm:4}
  Let $A$ be a $\modZ $-module, and let $S'\subset S\subset \modN $ be subsets.  Suppose
  that for each $n\in S$ there is a sequence $S'\ni n'\Leftrightarrow _A\cdots\Leftrightarrow _An$ in
  $S$.  Then the homomorphism
  $\rho ^A_{S,S'}\zzzcolon \AqS\rightarrow \Aq{S'}$ induced by $\id_{A[q]}$ is injective.
\end{theorem}

\begin{proof}
  One way to prove Theorem~\ref{thm:4} is to modify
  Section~\ref{sec:lemma-compl-polyn} and the proof of
  Theorem~\ref{thm:17}.  We roughly sketch the necessary
  modifications.  Section~\ref{sec:lemma-compl-polyn} is generalized
  as follows.  For two elements $f,g\in \modM _R$ and an $R$-module, we
  write $f\Rightarrow _Ag$ if $f\overset{I}\Rightarrow _Ag$ for some ideal $I$ with $A$
  being $I$-adically separated.  Then Proposition~\ref{thm:20} with
  $R$ replaced with an $R$-module $A$ holds.  Generalizations of
  Theorem~\ref{thm:21} and Corollary~\ref{thm:14} to $R$-modules is
  straightforward.  Theorem~\ref{thm:4} follows immediately from the
  generalized version of Corollary~\ref{thm:21}.

  Alternatively, we can use Theorem~\ref{thm:17} as follows.
  Since the case $A=0$ is trivial, we assume not.
  Let
  $A'=\modZ \oplus A$ be the ring with the multiplication
  $(m,a)(n,b)=(mn,mb+na)$ and with the unit $(1,0)$.  Then for
  $m,n\in \modN $ we have $m\Leftrightarrow _An$ if and only if $m\Leftrightarrow _{A'}n$.  Hence we can
  apply Theorem~\ref{thm:17} to obtain the injectivity of $\rho ^{A'}_{S,S'}$.
  We can identify $\rho ^{A'}_{S,S'}$
  with the direct product
  \[
  \rho ^{\modZ }_{S,S'}\oplus\rho ^{A}_{S,S'}\zzzcolon \ZqS\oplus\AqS\rightarrow \Zq{S'}\oplus\Aq{S'}.
  \]
  Hence $\rho ^A_{S,S'}$ is injective.
\end{proof}

\subsection{Non-surjectivity of $\rho ^\modZ _{\modN ,\{n\}}$}
\label{sec:non-surjectivity}

\begin{proposition}
  \label{thm:10}
  \begin{enumerate}
  \item If $m,n\in \modN $, $m\Leftrightarrow _\modZ n$, and $m\neq n$, then the homomorphism
    $\rho ^\modZ _{\{m,n\},\{m\}}\zzzcolon \Zq{\{m,n\}}\rightarrow \Zq{\{m\}}$ is {\em not}
    surjective.
  \item If $m|n$ and $m\neq n$, then the homomorphism
    $\rho ^\modZ _{\langle n\rangle ,\langle m\rangle }\zzzcolon \Zq{\langle n\rangle }\rightarrow \Zq{\langle m\rangle }$ is {\em not} surjective.
  \item For each nonempty, finite subset $S\subset \modN $,  the homomorphism
    $\rho ^\modZ _{\modN ,S}\zzzcolon \ZqN\rightarrow \ZqS$ is {\em not} surjective.
  \end{enumerate}
\end{proposition}

\begin{proof}
  (1) We have $m/n=p^e$ for some
  prime $p$ and an integer $e\neq 0$.  Consider the following
  commutative diagram of natural homomorphisms.
  \[
  \begin{CD}
    \Zq{\{m,n\}} @>\rho ^\modZ _{\{m,n\},\{m\}}>> \Zq{\{m\}} \\
    @VVV @VVbV \\
    \modZ [q]/(\Phi _n(q)) @>>c> \modZ _p[q]/(\Phi _n(q))
  \end{CD}
  \]
  It follows from $\modZ _p[q]/(\Phi _n(q))\simeq
  \varprojlim_j{\modZ [q]/(\Phi _n(q), p^j)}$, $\Phi _m(q)\in \sqrt{(\Phi _n(q),p)}$,
  and $p\in (\Phi _m(q),\Phi _n(q))$ that $b$ is a well-defined, surjective
  homomorphism.  Since $c$ is not surjective, $\rho ^\modZ _{\{m,n\},\{m\}}$
  is not surjective.

  (2) We may assume that $n=pm$ for a prime $p$.  The case $m=1$ is
  contained in (1) above.  There are isomorphisms
  $\Zq{\langle m\rangle }\simeq \modZ [q^m]^{\langle 1\rangle }\otimes _{\modZ [q^m]}\modZ [q]$ and
  $\Zq{\langle pm\rangle }\simeq \modZ [q^m]^{\langle p\rangle }\otimes _{\modZ [q^m]}\modZ [q]$ induced by the
  isomorphism $\modZ [q]\simeq\modZ [q^m]\otimes _{\modZ [q^m]}\modZ [q]$.  Then the case
  $m=1$ implies the non-surjectivity of $\rho ^\modZ _{\langle pm\rangle ,\langle m\rangle }$.

  (3) This follows from (2) above, since $\rho ^\modZ _{\modN ,S}$ factors
      through $\rho ^\modZ _{\langle n\rangle ,\langle m\rangle }$ for some $m,n$ with $m|n$ and
      $m\neq n$.
\end{proof}

\subsection{The ring $\QqS$}
\label{sec:structure-qqh}

The structure of $\QqS$ for $S\subset \modN $ is quite contrasting to that of $\ZqS$.
Note that $\ZqS$ embeds into $\QqS$ by Lemma~\ref{lem:1}.  (The
following remarks holds if we replace $\modQ $ with any ring $R$ such that each
element of $S$ is a unit in $R$.)

Note that if $m,n\in S$, $m\neq n$, then $(\Phi _m(q)^i,\Phi _n(q)^j)=(1)$ in
$\modQ [q]$ for any $i,j\ge 0$.  Consequently, for each
$f(q)=\prod_{n\in S}\Phi _n(q)^{\lambda (n)}\in \PHI_S^*$ with $\lambda (n)\ge 0$
we have by the Chinese Remainder Theorem
\[
\modQ [q]/(f(q)) \simeq \prod_{n\in S}\modQ [q]/(\Phi _n(q)^{\lambda (n)}).
\]
Taking the inverse limit, we obtain an isomorphism
\[
  \QqS\simeqto \prod_{n\in S}\Qqn.
\]
Since each $\Qqn$ is not zero, it follows that $\QqS$ is not an
integral domain if $|S|>1$.  It also follows
that $\rho ^\modQ _{S,S'}\zzzcolon \QqS\rightarrow \Qq{S'}$ is not
injective (but surjective) for each $S'\subsetneq S$.
Since for each $n\in S$ the
(surjective) homomorphism $\Qqn\rightarrow \modQ [q]/(\Phi _n(q))$ is not injective,
the homomorphism $\tau ^\modQ _{S,S}\zzzcolon \QqS\rightarrow \modQ [q]/(\Phi _n(q))$ is not
injective.


\begin{thebibliography}{99}
\bibitem{Habiro:rims2001} K.~Habiro, {\it On the quantum $sl_2$
    invariant of knots and integral homology spheres}, preprint.

\bibitem{Habiro:in-preparation} \bysame, in preparation.

\bibitem{Lawrence} R.~J.~Lawrence, {\it Asymptotic expansions of
    Witten-Reshetikhin-Turaev invariants for some simple
    $3$-manifolds}, J. Math. Phys. {\bf 36} (1995), no.~11,
  6106--6129.

\bibitem{Lawrence-Zagier} \bysame and D.~Zagier, {\it Modular
    forms and quantum invariants of $3$-manifolds}, Asian J. Math.
  {\bf 3} (1999), no.~1, 93--107.

\bibitem{Le} T.~T.~Q.~Le, {\it Quantum invariants of $3$-manifolds:
    integrality, splitting, and perturbative expansion}, preprint,
  {\tt math.QA/0004099}.
  
%
%
%

\bibitem{Murakami} H.~Murakami, {\it Quantum ${\rm
      SU}(2)$-invariants dominate Casson's ${\rm SU}(2)$-invariant},
  Math. Proc. Cambridge Philos. Soc. {\bf 115} (1994), no.~2,
  253--281.

\bibitem{Ohtsuki} T.~Ohtsuki, {\it A polynomial invariant of
    integral homology $3$-spheres}, Math. Proc. Cambridge Philos.
  Soc. {\bf 117} (1995), no.~1, 83--112.

\bibitem{Ohtsuki:problemlist} \bysame (ed.), {\it Problems on invariants of
    knots and $3$-manifolds}, preprint, available at {\tt
    http://www.ms.u-tokyo.ac.jp/\~{}tomotada/proj01}.

\bibitem{Reshetikhin-Turaev} N.~Reshetikhin and V.~G.~Turaev, {\it
    Invariants of $3$-manifolds via link polynomials and quantum
    groups}, Invent. Math.  {\bf 103} (1991), no.~3, 547--597.

\bibitem{Rozansky} L.~Rozansky, {\it On $p$-adic properties
    of the Witten-Reshetikhin-Turaev invariant}, preprint, {\tt
    math.QA/9806075}.

\bibitem{Sehgal} S.~K.~Sehgal, ``Units in integral group rings'',
  Pitman Monogr. Surveys Pure Appl. Math., Longman, Essex, 1993.

\bibitem{Witten} E.~Witten, {\it Quantum field theory and the Jones
    polynomial}, Comm. Math. Phys.  {\bf 121}  (1989),  no.~3, 351--399.

\bibitem{Zagier} D.~Zagier, {\it Vassiliev invariants and a strange
    identity related to the Dedekind eta-function}, Topology {\bf 40}
  (2001), no.~5, 945--960.
\end{thebibliography}
\end{document}